\documentclass[12pt]{article}%
\usepackage{amssymb}
\usepackage{amsfonts}
\usepackage{amsmath}
\usepackage{graphicx}%
\begin{document}

\begin{center}
{\large INTERVAL\ DECOMPOSITION LATTICES }

{\large ARE BALANCED}

\medskip

\textsc{Stephan Foldes}

{\footnotesize \textit{Institute of Mathematics, Tampere University of
Technology}}

{\footnotesize \textit{33101 Tampere, Finland}}

{\small e-mail:}{\footnotesize \textit{\ stephan.foldes@tut.fi}}

\smallskip

\textsc{and}

\smallskip

\textsc{S\'{a}ndor Radeleczki}

{\footnotesize \textit{Institute of Mathematics, University of Miskolc}}

{\footnotesize \textit{3515 Miskolc-Egyetemv\'{a}ros, Hungary}}

{\small e-mail:}{\footnotesize \textit{\ matradi@uni-miskolc.hu}}
\end{center}

\smallskip

\begin{center}
\textbf{Abstract}\medskip
\end{center}

{\footnotesize Intervals in binary or n-ary} {\footnotesize relations or other
discrete structures generalize the concept of interval in a linearly ordered
set. Join-irreducible partitions into intervals are characterized in the
lattice of all interval decompositions of a set, in a general sense of
intervals defined axiomatically. This characterization is used to show that
the lattice of interval decompositions is balanced.}

\medskip


\medskip

\noindent\textbf{2000 Mathematics Subject Classification: }%
{\footnotesize Primary 06B05, 06A15, 06C10; Secondary 05C99, 03C99.}

\bigskip

\begin{center}
\textsc{1. Preliminaries}
\end{center}

\medskip

Decompositions into intervals were first studied by Hausdorff [12,13], in the
context of linearly ordered sets, then extended to partially ordered sets,
graphs (see Sabidussi [19] ), appearing in particular in the study of
comparability graphs [9] ). The concept of decomposition was extended to
hypergraphs and directed graphs by D\"{o}rfler and Imrich [4] and D\"{o}rfler
[3], and to higher arity relational structures by Fraiss\'{e} [6,7,8]. A
general, abstract theory of decompositions was first presented by M\"{o}hring
and Radermacher [17] and M\"{o}hring [16]. Under mild stipulations about what
is to be considered an interval, interval decompositions constitute a complete
lattice. In [17] it was proved that this lattice is semimodular whenever it is
of a finite length. This result was extended in [5] to arbitrary interval
decomposition lattices, and their meet-irreducible elements were also
described. Another proof for the semimodularity can be found in [14]. We note
that in graph theory, interval decompositions are closely related to the
transitive orientation problem (see [9] and [15]). In the present paper, we
prove further properties of the lattice of interval decompositions. First, we
characterise the join-irreducible elements in this lattice, and using this
result we show that the lattice is balanced. As a consequence, several other
properties of the lattice of interval decompositions are deduced. \medskip

A closure system $(V,$ $\mathcal{Q})$, $\mathcal{Q}\subseteq\mathcal{P}(V)$ is
called \textit{algebraic} if the union of any chain of closed sets is closed.
An\textit{\ interval system} $(V,\mathcal{I})$ was defined in [5] as an
algebraic closure system with the following properties:\medskip

\noindent(I$_{0}$) $\{x\}\in\mathcal{I}$ for all $x\in V$ and $\varnothing
\in\mathcal{I}$,

\noindent(I$_{1}$) $A,B\in\mathcal{I}$ and $A\cap B\neq$ $\varnothing$ imply
$A\cup B\in\mathcal{I}$,

\noindent(I$_{2}$) For any $A,B\in\mathcal{I}$ the relations $A\cap B\neq$
$\varnothing$, $A\nsubseteqq B$ and $B\nsubseteqq A$ imply $A\setminus
B\in\mathcal{I}$ (and $B\setminus A\in\mathcal{I}$).

\medskip

\noindent Examples of interval systems given in (5) include modules of graphs
and relational intervals. These latter include order intervals in linearly
ordered sets.

\medskip

A set $A\in\mathcal{Q}$ is called a\textit{\ strong set }in the closure system
$(V,\mathcal{Q})$, if for any $B\in\mathcal{Q}$, $A\cap B\neq\varnothing$
implies $A\subseteqq B$ or $B\subseteqq A$. Let $\mathcal{S}$ stand for the
set of the all strong sets in $(V,\mathcal{Q})$; then $(V,\mathcal{Q})$ is an
algebraic closure system satisfying conditions (I$_{1}$) and (I$_{2}$). Let
$(V,\mathcal{Q})$ satisfy condition (I$_{0}$). Then $\varnothing$ and any
singleton $\{a\}$, $a\in V$ are strong sets, and hence $(V,\mathcal{S})$ is an
interval system. $\varnothing$, $V$ and the singletons $\{a\}$, $a\in V$ are
called \textit{improper} \textit{strong sets.}

\medskip

We note that restricting a closure system $(V,\mathcal{Q})$ to a nonempty set
$A\subseteq V$ we obtain again a closure system $(A,\mathcal{Q}_{A})$ with
$\mathcal{Q}_{A}=\{Q\cap A\mid Q\in\mathcal{Q\}}$. Clearly, for any
$A\in\mathcal{Q}$ we have $\mathcal{Q}_{A}\subseteq\mathcal{Q}$, and
$(A,\mathcal{Q}_{A})$ is an interval system whenever $(V,\mathcal{Q})$ is an
interval system.

\bigskip

\noindent\textbf{Definition 1.1. }A \textit{decomposition} in a closure system
$(V,\mathcal{Q})$ is a partition $\pi=\{A_{i}\mid i\in I\}$ of the set $V$
such that $A_{i}\in\mathcal{Q}$, for all $i\in I.$ If $(V,\mathcal{Q})$ is an
interval system, then $\pi$ is called an \textit{interval decomposition}. The
set of all decompositions in $(V,\mathcal{Q})$ is denoted by $\mathcal{D}%
(V,\mathcal{Q})$.

\bigskip

\noindent Let Part$(V)$ denote the lattice of all partitions of $V$. Since
$\mathcal{D}(V,\mathcal{Q})\subseteq\ $Part$(V)$, it is ordered by refinement,
where for any $\pi_{1},\pi_{2}\in\mathcal{D}(V,\mathcal{Q})$, $\pi_{1}\leq
\pi_{2}$ holds if and only if any block of $\pi_{2}$ is the union of some
blocks of $\pi_{1}$. Moreover, in [5] we proved the following.

\bigskip

\noindent\textbf{Proposition 1.2. }\emph{Let} $(V,\mathcal{Q})$ \emph{be a
closure system. Then} $\mathcal{D}(V,\mathcal{Q})$ \emph{is a complete lattice
with the greatest element }$\nabla=\{V\}$. \emph{If }$(V,\mathcal{Q})$
\emph{is algebraic and} \emph{satisfies condition }(I$_{0}$)\emph{, then}
$\mathcal{D}(V,\mathcal{Q})$ \emph{is a complete sublattice of }Part$(V)$
\emph{if and only if it satisfies condition }(I$_{1}$).

\bigskip

\noindent\textbf{Example 1.3. } If $T=(V,E)$ is a finite tree, then the vertex
sets of its subtrees form a closure system $(V,\mathcal{Q})$ which satisfies
conditions (I$_{0}$) and (I$_{1}$). Then $\mathcal{D}(V,\mathcal{Q})$ is a
finite sublattice of Part$(V)$, according to Proposition 1.2. We prove that
$\mathcal{D}(V,\mathcal{Q)}$ is a Boolean lattice isomorphic to $(\mathcal{P}%
(E),\subseteq)$.

Indeed, given $S\subseteq E$ define $\pi(S)$ as the equivalence relation on
$V$ in which two vertices are equivalent if they are connected in the tree $T
$ by a path containing only edges from $S$. Since the classes of $\pi(S)$
induce subtrees of $T$, $\pi(S)$ is a decomposition in $(V,\mathcal{Q)}$. Then
the isomorphism of $\mathcal{P}(E)$ to the lattice $\mathcal{D}(V,\mathcal{Q)}%
$ is given by the mapping
\[
S\longmapsto\pi(S)\text{.}%
\]

\noindent The following result was proved in [5]:

\bigskip

\noindent\textbf{Proposition 1.4. }\emph{If }$(V,\mathcal{Q})$ \emph{be an
algebraic closure system satisfying condition }(I$_{1}$)\emph{, then}
$\mathcal{D}(V,\mathcal{Q})$ \emph{is an algebraic semimodular lattice.}

\medskip

Therefore, for an\textit{\ }interval system $(V,\mathcal{I})$, the lattice
$\mathcal{D}(V,\mathcal{I})$ is always an algebraic semimodular sublattice of
Part$(V)$.\medskip

\noindent\textbf{Remark 1.5.} Let $(V,\mathcal{Q})$ be a closure system
satisfying (I$_{0}$). Then clearly $\triangle=\{\{x\}\mid x\in V\}$ is the
least element of $\mathcal{D}(V,\mathcal{Q})$, and to any $A\in\mathcal{Q}%
\setminus\{\varnothing\}$ corresponds the decomposition%
\[
\pi_{A}=\{A\}\cup\{\{x\}\mid x\in V\setminus A\}\text{.}%
\]

\noindent Moreover, if $\pi=\{A_{i}\mid i\in I\}\in\mathcal{D}(V,\mathcal{Q}%
)$, then%
\begin{equation}
\pi=%
{\textstyle\bigvee}
\{\pi_{A_{i}}\mid i\in I\} \tag{1}%
\end{equation}

\noindent where $%
{\textstyle\bigvee}
$ means the join in the complete lattice $\mathcal{D}(V,\mathcal{Q})$.

\medskip

A decomposition $\pi=\{A_{i}\mid i\in I\}$ in a closure system $(V,\mathcal{Q}%
)$ is called a \textit{strong decomposition} if every $A_{i}$, $i\in I$ is a
strong set in $(V,\mathcal{Q})$. Since the strong decompositions in
$(V,\mathcal{Q})$ can be considered also as decompositions in the closure
system $(V,\mathcal{S})$, they form a complete lattice $\mathcal{D}%
(V,\mathcal{S})$ whose greatest element is $\triangledown=\{V\}$.

An element $a$ of a lattice $L$ is called \textit{standard }(see Gr\"{a}tzer
[10]), if%
\[
x\wedge(a\vee y)=(x\wedge a)\vee(x\wedge y)\text{ holds for all }x,y\in
L\text{.}%
\]

\noindent The standard elements of $L$ form a distributive sublattice of $L$
denoted by S$(L)$. The following result was proved also in [5]:

\bigskip

\noindent\textbf{Theorem 1.6. }Let $(V,\mathcal{Q})$ \emph{be a closure
system. Then the} \emph{strong decompositions in }$(V,\mathcal{Q})$ \emph{are
standard elements of} $\mathcal{D}(V,\mathcal{Q})$ \emph{and} $\mathcal{D}%
(V,\mathcal{S})$ \emph{is a distributive sublattice of} $\mathcal{D}%
(V,\mathcal{Q})$ \emph{and of} Part$(V)$.

\bigskip

A set $A\in\mathcal{Q}$ of a closure system $(V,\mathcal{Q})$ is called
\textit{fragile} if it is the union of two disjoint nonempty members of
$\mathcal{Q}$, otherwise $A$ is called \textit{nonfragile}. This generalizes
the concept of fragility studied by Habib and Maurer [11] in the context of
the module systems of graphs. In view of [5], if $(V,\mathcal{Q})$ is an
interval system, then any nonfragile interval\emph{\ }$A\in\mathcal{Q}$ is a
strong set.

\bigskip

\begin{center}
\textsc{2. Completely join-irreducible elements in }$\mathcal{D}%
(V,\mathcal{I})$
\end{center}

\bigskip

An element $j\in L\setminus\{0\}$ of a complete lattice $L$ is
called\textit{\ completely join-irreducible }if for any system of elements
$x_{i}\in L$, $i\in I$ the equality $j=\bigvee\{x_{i}\mid i\in I\}$ implies
$p=x_{k}$ for some $k\in I$. Let $J(L)$ stand for the set of completely
join-irreducible elements of $L$. The completely meet-irreducible elements of
$L$ are defined dually, and their set is denoted by $M(L)$. Let us define
$a_{\ast}:=%
{\textstyle\bigvee}
\{x\in L\mid x<a\}$ for any $a\in L\setminus\{0\}$, and $a^{\ast}:=$ $%
{\textstyle\bigwedge}
\{x\in L\mid x>a\}$, for any $a\in L\setminus\{1\}$. Denoting by $\prec$ the
covering relation in a lattice $L$, we can observe that
\[
j\in L\setminus\{0\}\text{ is completely join-irreducible}\Leftrightarrow
j_{\ast}<j\Leftrightarrow j_{\ast}\prec j\text{, and}%
\]%
\[
m\in L\setminus\{1\}\text{ is completely meet-irreducible}\Leftrightarrow
m<m^{\ast}\Leftrightarrow m\prec m^{\ast}\text{.}%
\]

\medskip

In this section we characterise the completely join-irreducible elements in
the lattice $\mathcal{D}(V,\mathcal{I)}$ of interval decomposition, and we
show that they are closely related to the strong sets of $(V,\mathcal{I)}$.

\medskip

The following result of [5] will be useful in our proofs.

\bigskip

\noindent\textbf{Lemma 2.1. }\emph{Let }$\pi_{1}=\{B_{j}\mid j\in
J\}$\emph{\ and }$\pi_{2}=\{A_{i}\mid i\in I\}$\emph{\ be two decompositions
in a closure system} $(V,\mathcal{Q})$. \emph{If }$\pi_{1}\prec\pi_{2}$
\emph{holds in} $\mathcal{D}(V,\mathcal{Q})$\emph{\ then there exists a unique
}$k\in I$\emph{\ and }$J_{k}\subseteqq J$\emph{\ with at least two elements
such that }$A_{k}=\ \underset{j\in J_{k}}{\bigcup}B_{j}$\emph{\ and }$A_{i}%
\in\pi_{1}$\emph{\ for all }$i\in I\setminus\{k\}$\emph{.}

\medskip

Now, we are prepared to prove the following result

\bigskip

\noindent\textbf{Lemma 2.2. }\emph{Let} $(V,\mathcal{Q})$ \emph{be a closure
system satisfying condition} (I$_{0}$) \emph{and} $\pi\in\mathcal{D}%
(V,\mathcal{Q})$. \emph{Then }$\pi$\emph{\ is completely join-irreducible in
}$\mathcal{D}(V,\mathcal{Q})$ \emph{if and only if there exists a nonempty
closed set} $A\in\mathcal{Q}$ \emph{with }$\pi=\pi_{A}$\emph{\ and such that
}$A$\emph{\ admits a greatest proper decomposition into closed sets.}

\bigskip

\noindent\textit{Proof.} Assume that $\pi=\{A_{i}\mid i\in I\}$ is a
completely join-irreducible element in $\mathcal{D}(V,\mathcal{Q})$. Since
$\pi=%
{\textstyle\bigvee}
\{\pi_{A_{i}}\mid i\in I\}$ according to (1), we obtain that $\pi=\pi_{A_{k}}%
$, for some $k\in K$. Then, in view of Lemma 2.1, $\pi_{\ast}=\{B_{j}\mid j\in
J\}\cup\{\{x\}\mid x\in V\setminus A_{k}\}$, where $B_{j}\in\mathcal{Q}$,
$\mid J\mid\geq2$, and $A_{k}$ equals to the disjoint union $\underset{j\in
J}{\bigcup}B_{j}$. Clearly, $\mu=\{B_{j}\mid j\in J\}$ is a proper
decomposition of $(A,\mathcal{Q}_{A_{k}})$. Let $\nu=\{C_{t}\mid t\in T\}$ be
an arbitrary decomposition of $(A,\mathcal{Q}_{A_{k}})$ such that $\nu
\neq\{A_{k}\}$. Then it is easy to see that $\nu^{+}=\{C_{t}\mid t\in
T\}\cup\{\{x\}\mid x\in V\setminus A_{k}\}$ is a decomposition in
$\mathcal{D}(V,\mathcal{Q})$ and $\nu^{+}<\pi_{A_{k}}=\pi$. Since $\pi$ is
completely join-irreducible, we get $\nu^{+}<$ $\pi_{\ast}$. Hence the
partition $\nu=\{C_{t}\mid t\in T\}$ of $A$ is a refinement of the partition
$\mu=\{B_{j}\mid j\in J\}$ of $A$. Thus $\nu\leq\mu$ holds in $\mathcal{D}%
(A,\mathcal{Q}_{A_{k}})$, and this means that $A_{k}$\emph{\ }admits $\mu$ as
a greatest proper decomposition.

Conversely, let $A\in\mathcal{Q}\setminus\{\emptyset\}$ be a closed set that
admits a greatest proper decomposition $A=\{B_{j}\mid j\in J\}$, $\mid
J\mid\geq2$. We prove that $\pi_{A}$ is a completely join-irreducible element
in $\mathcal{D}(V,\mathcal{Q})$. Since $B_{j}\in\mathcal{Q}_{A}\subseteq
\mathcal{Q}$ and $B_{j}\neq\emptyset$, for each $j\in J$, $\{B_{j}\mid j\in
J\}\cup\{\{x\}\mid x\in V\setminus A\}$ is a decomposition in $\mathcal{D}%
(V,\mathcal{Q})$.

Now, let $\mu\in\mathcal{D}(V,\mathcal{Q})$, $\mu<\pi_{A}$ arbitrary. Since
the partition $\mu$ is a refinement of $\pi_{A}$, it is has the form
$\mu=\{C_{t}\mid t\in T\}\cup\{\{x\}\mid x\in V\setminus A\}$, where $C_{t}%
\in\mathcal{Q}$, $\mid T\mid\geq2$, and $A=\ \underset{t\in T}{\bigcup}C_{t}$.
Then $\{C_{t}\mid t\in T\}$ is a proper decomposition in $(A,\mathcal{Q}_{A}%
)$, and hence $\{C_{t}\mid t\in T\}\leq\{B_{j}\mid j\in J\}$. Because this
result yields $\mu\leq\{B_{j}\mid j\in J\}\cup\{\{x\}\mid x\in V\setminus
A\}$, we deduce $(\pi_{A})_{\ast}\leq\{B_{j}\mid j\in J\}\cup\{\{x\}\mid x\in
V\setminus A\}<\pi_{A}$ , and this implies that $\pi_{A}$ is completely
join-irreducible.\hfill$\blacksquare$

\bigskip

\noindent\textbf{Proposition 2.3. }\emph{Let} $(V,\mathcal{I})$ \emph{be an
interval system. Then }$\pi$\emph{\ is a completely join-irreducible element
in }$\mathcal{D}(V,\mathcal{I})$ \emph{if and only if }$\pi=\pi_{A} $\emph{,}
\emph{where} $A$ \emph{is an} \emph{interval} \emph{admitting a greatest
proper decomposition} $\{B_{j}\mid j\in J\}$, $\mid J\mid\geq2$ \emph{such
that each} $B_{j}$, $j\in J$ \emph{is a strong set in }$(V,\mathcal{I}%
)$\emph{. If the later condition holds with }$\mid J\mid\geq3$, \emph{then}
$A$ \emph{is strong.}

\bigskip

\noindent\textit{Proof. }Let\textit{\ }$A\in\mathcal{I}\setminus\{\emptyset\}$
be an interval such that $\mu=\{B_{j}\mid j\in J\}$, $\mid J\mid\geq2$ is a
greatest proper decomposition of it. Then $B_{j}\in(V,\mathcal{I})$ for all
$j\in J$. If $\mid J\mid\geq3$, then $A$ can not be the union of two disjoint
nonempty members of $\mathcal{I}$. Therefore, $A$ is nonfragile and hence it
is a strong set according to [5]. In view of Lemma 2.2, to prove our statement
it is enough to show that each $B_{j}$, $j\in J$ is strong.

First, assume that $\mid J\mid\geq3$. Then $A$ is a strong set. Let
$C\in\mathcal{I}$ such that $C\cap B_{j}\neq\emptyset$, for some $j\in J$.
Since $A$ is strong, now $C\cap A\neq\emptyset$ implies that either
$A\subseteq C$ or $C\varsubsetneqq A$ holds. In the first case $B_{j}\subseteq
C$. If $C\varsubsetneqq A$, then $\nu=\{C\}\cup\{\{x\}\mid x\in A\setminus
C\}$ is a decomposition in $(A,\mathcal{I}_{A})$ and $\nu\neq\{A\}$. Hence
$\nu\leq\mu$, and this implies $C\subseteq B_{j}$, because $C\cap B_{j}%
\neq\emptyset$. Therefore, $B_{j}$ is a strong set of $(V,\mathcal{I})$.

Let $\mid J\mid=2$. Then $A=B_{1}\cup B_{2}$ and $\mu=\{B_{1},B_{2}\}$ is the
greatest proper decomposition in $(A,\mathcal{I}_{A})$. Suppose that $C\cap
B_{1}\neq\emptyset$ and $B_{1}\nsubseteqq C$. Then also $C\cap A\neq\emptyset$
and $A\nsubseteqq C$. Since $C\cap A$, $A\setminus C\in\mathcal{I}%
\setminus\{\emptyset\}$, we get that $\rho=\{C\cap A$, $A\setminus C\}$ is
proper decomposition in $(A,\mathcal{I}_{A})$. Hence $\rho\leq\mu$. Since
$\rho$ a maximal proper partition of the set $A$, we obtain $\rho=\mu$. Then
$B_{1}=C\cap A$, as otherwise $B_{1}=$ $A\setminus C$ would imply $C\cap
B_{1}=\emptyset$, a contradiction. Thus $B_{1}\subseteq C$ and this means that
$B_{1}$ is a strong set. The fact that $B_{2}$ is strong is proved
similarly.\hfill$\blacksquare$

\bigskip

As an immediate consequence of Proposition 2.3 we obtain:

\bigskip

\noindent\textbf{Corollary 2.4. }\emph{Let} $(V,\mathcal{I})$ \emph{be an
interval system. If }$\pi$\emph{\ is a completely join-irreducible element in
}$\mathcal{D}(V,\mathcal{I})$\emph{, then }$\pi_{\ast}$\emph{\ is a strong
decomposition.}

\bigskip

A lattice $L$ is called \textit{geometric}, if it is atomistic, semimodular
and algebraic. It is well-known that the geometric lattices are also dually atomistic.

\bigskip

\noindent\textbf{Corollary 2.5.} \emph{Let} $(V,\mathcal{I})$ \emph{be an
interval system such that }$\mathcal{D}(V,\mathcal{I})$ \emph{is of finite
length. Then the following assertions are equivalent:}

\noindent(i) $\mathcal{D}(V,\mathcal{I})$ \emph{is an atomistic lattice;}

\noindent(ii) $\mathcal{D}(V,\mathcal{I})$ \emph{is a geometric lattice;}

\noindent(iii) $\mathcal{D}(V,\mathcal{I})$ \emph{is a dually atomistic
lattice;}

\noindent(iv) $\mathcal{D}(V,\mathcal{I})$ \emph{has no proper strong
intervals.}

\bigskip

\noindent\textit{Proof. }(i)$\Rightarrow$(ii) is clear, because $\mathcal{D}%
(V,\mathcal{I})$ is an algebraic semimodular lattice, according to Proposition
1.5. The implication (ii)$\Rightarrow$(iii) is obvious, and (iii)$\Rightarrow
$(iv) follows from [5, Corollary 3.7].

(iv)$\Rightarrow$(i). Since $\mathcal{D}(V,\mathcal{I})$ is a lattice of
finite length, any element of it is a join of some completely join-irreducible
elements (see e.g. [1]). Now, assume that $(V,\mathcal{I})$ has no proper
strong intervals, and let $\pi$ be a completely join-irreducible element of
$\mathcal{D}(V,\mathcal{I})$. Since $\pi>\triangle$, we get $\pi_{\ast}%
\geq\triangle$. Because by Corollary 2.4 $\pi_{\ast}\neq\{V\}$ is a strong
decomposition, in view of Proposition 2.3 we get that any block of $\pi_{\ast
}$ is of the form $\{a\}$, $a\in V$. Then $\pi_{\ast}=\triangle$. Since
$\triangle$ is the $0$-element of $\mathcal{D}(V,\mathcal{I})$, it follows
that $\pi$ is an atom. Hence $\mathcal{D}(V,\mathcal{I})$ is an atomistic
lattice.\hfill$\blacksquare$

\bigskip

\begin{center}
\textsc{3. Further properties of the lattice }$\mathcal{D}(V,\mathcal{I})$
\end{center}

\bigskip

Let $L$ be a lattice of finite length. $L$ is called a \textit{strong lattice}
if for any join-irreducible element $j\in J(L)$ and for all $x\in L$%
\[
j\leq j_{\ast}\vee x\text{ implies }j\leq x\text{.}%
\]

\noindent It is easy to see, that any atomistic lattice is strong. We say that
the lattice $L$ is \textit{dually strong}, if its dual $L^{(d)}$ is strong.
$L$ is called a \textit{consistent lattice}, if for any $j\in J(L)$ and each
$x\in L$, the element $x\vee j$ is a join-irreducible in the interval $[x,1]$.
If for any $j\in J(L)$ and $m\in M(L)$ with $j\nleq m$
\[
j\vee m=m^{\ast}\Leftrightarrow j\wedge m=j_{\ast}\text{.}%
\]

\noindent holds true, then $L$ is called a \textit{balanced} lattice. We say
that $L$ \textit{satisfies the Kurosh-Ore replacement property for
join-decompositions} ($\vee$-KORP, for short), if for every $a\in L$, and any
two\ irredundant join-decompositions
\[
a=j_{1}\vee...\vee j_{m}\text{ and }a=k_{1}\vee...\vee k_{n}\text{,}%
\]

\noindent with $j_{1},...,j_{m},k_{1},...,k_{n}\in J(L)$, each $j_{i}$ can be
replaced by a $k_{p}$ such that
\[
a=j_{1}\vee...\vee j_{i-1}\vee k_{p}\vee j_{i+1}.\vee...\vee j_{m}\text{.}%
\]

\medskip

\noindent\textbf{Remark 3.1. }It is well-known that any semimodular lattice of
finite length is dually strong (see e.g. Stern [20]). It belongs to the
folklore that a lattice $L$ of finite length is balanced if and only if both
$L$ and $L^{(d)}$ are strong. Crawley showed [1] that $L$ satisfies the $\vee
$-KORP if and only if $L$ is consistent. Let $L$ be of finite length. As it is
noted in [20], from the previous facts together with Walendziak result [21,
Thm.1] it follows the equivalence of the following assertions:\medskip

\noindent(a) $L$ \emph{is semimodular and has the }$\vee$-KORP.

\noindent(b) $L$ \emph{is semimodular and balanced;}

\noindent(c) $L$ \emph{is semimodular and consistent;}

\noindent(d) $L$ \emph{is semimodular and strong.}

\bigskip

\noindent\textbf{Theorem 3.2.} \emph{Let} $(V,\mathcal{I})$ \emph{be an
interval system. If the lattice }$\mathcal{D}(V,\mathcal{I})$ \emph{has finite
length, then it is a balanced lattice that has the }$\vee$-KORP.

\bigskip

\noindent\textit{Proof. }Since\textit{\ }$\mathcal{D}(V,\mathcal{I})$ is a
semimodular lattice of finite length, in order to prove our theorem, in view
of Remark 3.1, it suffices only to show that $\mathcal{D}(V,\mathcal{I})$ is
strong. Take any $j\in J(\mathcal{D}(V,\mathcal{I}))$ and $x\in\mathcal{D}%
(V,\mathcal{I})$ with $j\leq j_{\ast}\vee x$. Because any join-irreducible
element of lattice with finite length is also completely join-irreducible, in
view of Corollary 2.4, $j_{\ast}$ is a standard element in $\mathcal{D}%
(V,\mathcal{I})$. Thus we obtain.%
\begin{equation}
j=j\wedge(j_{\ast}\vee x)=(j\wedge j_{\ast})\vee(j\wedge x)=j_{\ast}%
\vee(j\wedge x)\text{.} \tag{2}%
\end{equation}

\noindent Since $j$ is join-irreducible and $j_{\ast}<j$, (2) implies
$j=j\wedge x$. Hence $j\leq x$, and this proves that $\mathcal{D}%
(V,\mathcal{I})$ is strong. \hfill$\blacksquare$

\bigskip

The above theorem has an additional consequence for finite interval
decomposition lattices.

\medskip

\noindent A \textit{tolerance} of a lattice $L$ is a reflexive and symmetric
relation $T\subseteq L^{2}$, compatible with the operations of $L$. A
\textit{block} of $T$ is a maximal set $B\subseteq L$ satisfying
$B^{2}\subseteq T$. Suppose that the lattice $L$ is of finite length. Then any
block $B$ of $T$ has the form of an interval $B=[u,v]$, $u,v\in L$, $u\leq v$,
and the compatibility property of $T$ makes it possible to build a "factor
lattice" $L/T$, whose elements are the blocks of $T$ (see Cz\'{e}dli [2]). $T$
is called a \textit{glued tolerance}, if it contains all covering pairs of
$L$. Since every intersection of glued tolerances of $L$ is again a glued
tolerance of $L$, there exists a least tolerance $\Sigma(L)$ comprising all
pairs $x\prec y$ in $L$, called the \textit{skeleton tolerance} of $L$. The
lattice $L$ is said to be \textit{glued by geometric lattices,} if all blocks
of $\Sigma(L)$ are geometric lattices.

Reuter [18] proved (see also [20; Thm. 4.6.8]) that for a finite lattice $L$,
the assertions (a),(b),(c) and (d) of Remark 3.1 are equivalent to the
following statement:

\medskip

\noindent(e) $L$ \emph{is glued by geometric lattices.}

\medskip

Therefore, by Theorem 3.2 we infer:

\medskip

\noindent\textbf{Corollary 3.3.} \emph{Let} $(V,\mathcal{I})$ \emph{be an
interval system such that }$\mathcal{D}(V,\mathcal{I})$ \emph{is a finite}
\emph{lattice. Then }$\mathcal{D}(V,\mathcal{I})$ \emph{is glued by geometric
lattices. }

\bigskip\bigskip

\begin{center}

\end{center}

Acknowledgements.

This work has been co-funded by Marie Curie Actions and supported by the
National Development Agency (NDA) of Hungary and the Hungarian Scientific
Research Fund (OTKA, contract number 84593), within a project hosted by the
University of Miskolc, Department of Analysis.

\bigskip

\includegraphics[height=15mm, width=20mm]{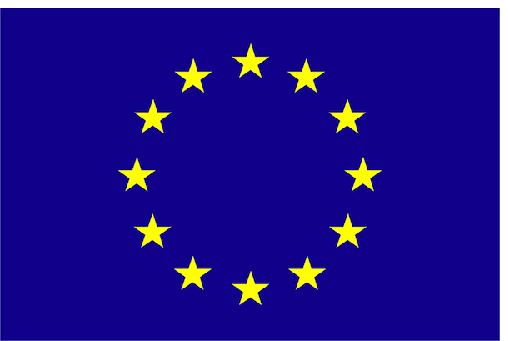}
\includegraphics[height=15mm, width=20mm]{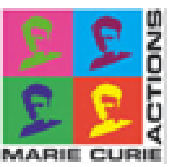}
\includegraphics[height=15mm, width=20mm]{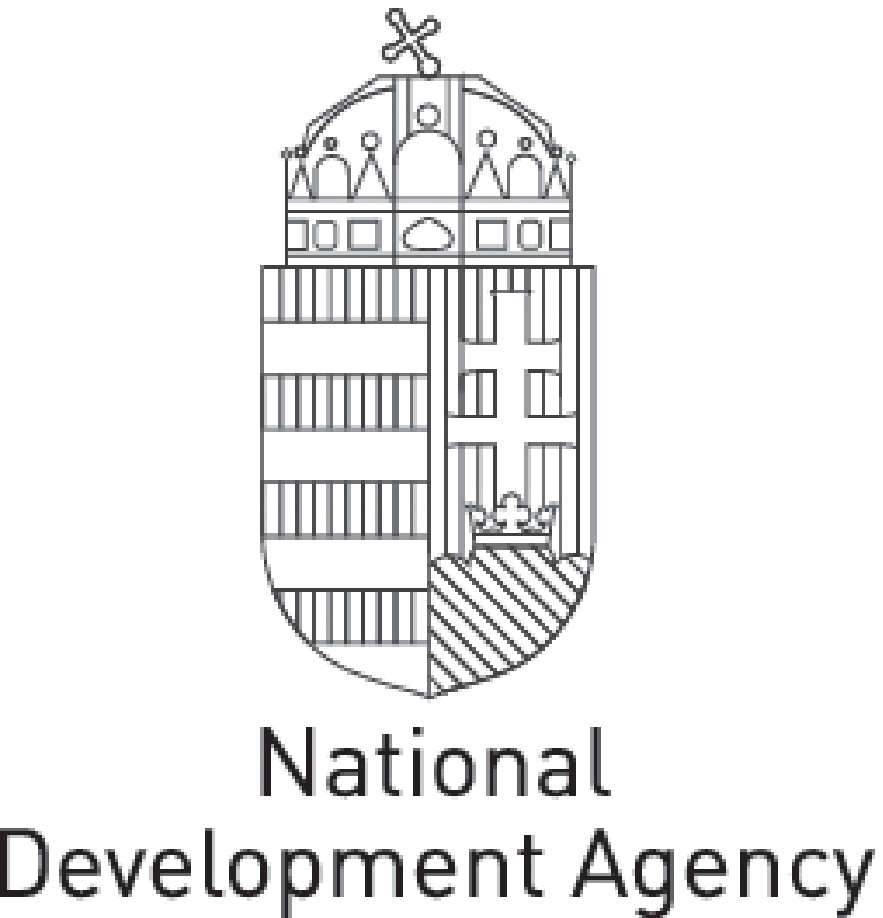} \includegraphics[height=15mm, width=20mm]{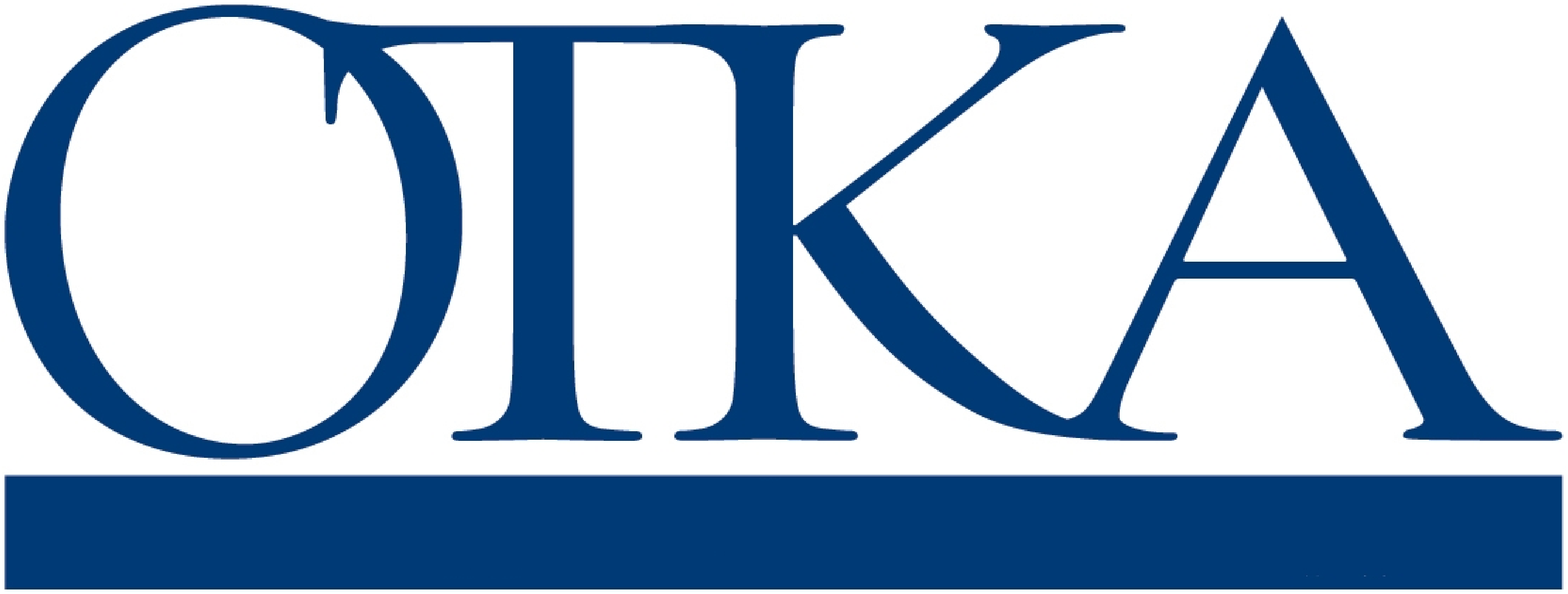}

\end{document}